\documentclass[11pt]{article}

       \font\tenmsb=msbm10
       \font\sevenmsb=msbm7
       \font\fivemsb=msbm5
       \catcode`\@=11
       \ifx\amstexloaded@\relax\catcode`\@=\active
       \endinput\else\let\amstexloaded@\relax\fi
       \def\spaces@{\space\space\space\space\space}
       \def\spaces@@{\spaces@\spaces@\spaces@\spaces@\spaces@}
       \def\space@.{\futurelet\space@\relax}
       \space@. %
       \def\Err@#1{\errhelp\defaulthelp@\errmessage{AmS-TeX error: #1}}
       \def\relaxnext@{\let\next\relax}
       \def\accentfam@{7}
       
       \def\noaccents@{\def\accentfam@{0}}
       \def\Cal{\relaxnext@\ifmmode\let\next\Cal@\else
       \def\next{\Err@{Use \string\Cal\space only in math mode}}\fi\next}
       \def\Cal@#1{{\Cal@@{#1}}}
       \def\Cal@@#1{\noaccents@\fam\tw@#1}

       \def\Bbb{\relaxnext@\ifmmode\let\next\Bbb@\else
       \def\next{\Err@{Use \string\Bbb\space only in math mode}}\fi\next}
       
       \def\Bbb@#1{{\Bbb@@{#1}}}
       \def\Bbb@@#1{\noaccents@\fam\msbfam#1}
       \newfam\msbfam
       \textfont\msbfam=\tenmsb
       \scriptfont\msbfam=\sevenmsb
       \scriptscriptfont\msbfam=\fivemsb
\usepackage{hyperref}
\usepackage{tikz}
\usepackage{pgflibraryarrows}
\usepackage{pgflibrarysnakes}
\usepackage{amsmath,amssymb}

\def\E{{\Bbb E}}
\def\N{{\Bbb N}}

\def\R{{\Bbb R}}

\def\C{{\Bbb C}}

\newtheorem{thm}{Theorem}
\newtheorem{thm*}{Theorem}
\newtheorem{lemma}{Lemma}[section]

\newtheorem{rmk}{Remark}[section]

\newtheorem{notation*}{Notation}

\@addtoreset{equation}{section}

\newcommand{\qed}{\nolinebreak\hfill\rule{2mm}{2mm}
\par\medbreak}

\newcommand{\beq}{\begin{equation} }
\newcommand{\eeq}{\end{equation} }

\newcommand{\calD}{\Cal D}

\begin{document}
\title{The exact Power Law for Buffon's needle landing near some Random Cantor Sets
}
\author{Shiwen Zhang\footnote
{Research of the author was supported in part by NSF grant DMS-1600065 .}}
\date{}
\maketitle
\begin{abstract}
 In this paper, we study the Favard length of some random Cantor sets of Hausdorff dimension 1.  We start with a unit disk in the plane and replace the unit disk by $4$ disjoint subdisks (with equal distance to each other) of radius $1/4$ inside and tangent to the unit disk. By repeating this operation in a self-similar manner and adding a random rotation in each step, we can generate a random Cantor set ${\cal D}(\omega)$. Let ${\cal D}_n$ be the $n$-th generation in the construction, which is comparable to the $4^{-n}$-neighborhood of ${\cal D}$. We are interested in the decay rate of the Favard length of
these sets ${\cal D}_n$ as $n\to\infty$, which is the likelihood (up to a constant) that ``Buffon's needle'' dropped randomly will fall into the $4^{-n}$-neighborhood of ${\cal D}$. It is well known in \cite{m} that the lower bound of the Favard length of ${\cal D}_n(\omega)$ is constant multiple of $n^{-1}$. We show that the upper bound of the Favard length of ${\cal D}_n(\omega)$ is $C n^{-1}$ for some $C>0$ in the average sense. We also prove the similar linear decay for the Favard length of ${\cal D}^d_n(\omega)$ which is the $d^{-n}$-neighborhood of a self-similar random Cantor set with degree $d$ greater than $4$. Notice in the non-random case where the self-similar set has degree greater than $4$, the best known result for the decay rate of the Favard length is $e^{-c\sqrt {\log n}}$.
\end{abstract}


\section{Introduction}\label{sec:intro}
Recall that the \emph{Favard length} of a planar
set $E\subset \C$ is defined by
\begin{equation}\label{Fav}
  {\rm Fav}(E)=\frac{1}{\pi}\int_0^\pi\big|{\rm Proj}_\theta\ E\big|{\rm d}\theta ,
\end{equation}
where ${\rm Proj}_\theta$ denotes the orthogonal projection onto the line having angle $\theta$ with the real axis and $|A|$ denotes the Lebesgue measure of a measurable set $A\subset \R$. For compact $E$, ${\rm Fav}(E)$ is the likelihood (up to a constant) that the ``Buffon's needle,'' a long line segment dropped with independent, uniformly distributed orientation and distance from the origin will intersect $E$. For this reason, ${\rm Fav}(E)$ is also called \emph{Buffon needle probability}.
A projection theorem of Besicovitch shows that if an unrectifiable set $E$ has finite 1-dimensional Hausdorff measure (${\cal H}^1(E)<\infty$), then the projection of $E$ to almost every
line through the origin has zero length, see \cite{m2}. Therefore, ${\rm Fav}(E)=0$ if $E$ is unrectifiable with finite length. The four-corner Cantor set ${\cal K}$ is one such example, which will be one of the central topics in this paper.

Let ${\cal C}$ be the middle-half cantor set. The four-corner Cantor set ${\cal K}$ is the Cartesian square ${\cal C}\times{\cal C}$. Let ${\cal C}_n$ be the n-th generation in the Cantor construction of ${\cal C}$ and let ${\cal K}_n={\cal C}_n\times{\cal C}_n$, then ${\cal K}=\bigcap{\cal K}_n$.  By Dominated Convergence Theorem, ${\rm Fav}({\cal K})=0$ implies $\lim_n{\rm Fav}({\cal K}_n)=0$. Due to the nested structure,  ${\rm Fav}({\cal K}_n)$ is comparable to the likelihood that
``Buffon's needle'' will land in a $4^{-n}$-neighborhood of ${\cal K}$. So people are interested in the decay rate of ${\rm Fav}({\cal K}_n)$ as $n\to\infty$. The lower bound ${\rm Fav}({\cal K}_n)\gtrsim n^{-1}$ follows from some general results relating the average projection length and Riesz 1-capacity by Mattila \cite{m}. The lower bound for ${\rm Fav}({\cal K}_n)$ was improved to be constant multiple of $n^{-1}\log n$ by Bateman and Volberg in \cite{bv}. The argument and result in \cite{bv} also hold true for Sierpi\'nski triangle and some general disk models, which are the best known results for these self-similar Cantor sets. Recently, some interesting examples of purely unrectifiable sets for which the Favard length decays at any rate chosen were considered by Wilson in \cite{w}. Wilson constructed non-self-similar Cantor sets, for which the Favard length of $\varepsilon$-neighborhoods decays arbitrarily with respect to $\varepsilon$.

 Back to the four-corner Cantor set, the first quantitative upper bound ${\rm Fav}({\cal K}_n)\le Ce^{-c\log^\ast n}$
 \footnote{$\log^\ast n$ denotes the
number of iterations of the $\log$ function needed to have $\log\cdots \log n \le 1$.} was obtained by Peres and Solomyak in \cite{ps} with extremely slow decay in $n$. More recently, Tao gave a quantitative statement of Besicovitch theorem \cite{t} which also provides a weaker upper bound for the Favard length in a more general setting.  But being rather general, the method in \cite{t}
does not give a good estimate for self-similar structures such as ${\cal K}_n$.
Vastly improved estimate with power decay  was proved by Nazarov, Peres and Volberg in \cite{npv}. They showed  ${\rm Fav}({\cal K}_n)\le C_p/n^{p} $ for any $p<\frac{1}{6}$ which is the best known result to our knowledge. The same type power estimate was proved for the $3^{-n}$-neighborhood of Sierpi{\'n}ski
gasket ${\cal S}_n$ by Bond-Volberg in \cite{bv1} and for rational product Cantor
sets with a tiling condition by {\L}aba-Zhai in \cite{lz}. More generalizations along these lines were considered in \cite{bv2,blv}. In \cite{bv2}, Bond and Volberg got the $e^{-c\sqrt {\log n}}$ upper bound, which holds true for some general self-similar Cantor sets with degree
 \footnote{The self-similar Cantor set is generated by using linear map interations. Degree here means the number of components contained in the image of the linear map.}
 greater than 4.  Bond, {\L}aba and Volberg in \cite{blv} studied the product rational cases and got the upper bound $n^{-c/ \log\log n}$. The estimates in \cite{bv2,blv}-which are worse than the expected polynomial estimate cast serious doubts on the general validity of the power estimate. The authors in \cite{blv} exposed the importance of  the structure of zeros of certain trigonometric polynomials for the decay of the Favard length. The question is related to understanding the zeros of fewnomials (polynomials with coefficients only 0 and 1) in algebra, see more discussion in \cite{blv}. This is  one important motivation to study the orthogonal projection of these particular Cantor sets.  The study of the decay rate of the Favard length and these trigonometric polynomials is also closely related to
the ``visibility'' of a planar set and Bernoulli convolutions. We refer readers to \cite{ps1,blz} for more details on these topics. Thorough reviews on orthogonal projections and more related topics about Favard length can also be found in \cite{pss,m3,l} and references therein.

Another important example which essentially inspires our paper is a random analogue of ${\cal K}_n$ considered in \cite{ps}. Recall that an alternative way to construct ${\cal K}_n$ is to replace the unit square by four subsquares of side length $1/4$ at its corners, and iterate this operation in a self-similar manner in each subsquare. In \cite{ps}, Peres and Solomyak considered the following random construction. Partition the unit square into four subsquares of side $1/2$. In each of these four subsquares, partition dyadically again and then choose a dyadic subsquare of side 1/4, uniformly at random. The union of these four squares is denoted as ${\cal R}_1$. Repeat this operation in each of these quarter subsquares randomly in a self-similar manner. At the $n$-th  step, it is a union of $4^n$ subsquares of side $4^{-n}$, denoted by ${\cal R}_n$. ${\cal R}=\bigcap_n{\cal R}_n$ will be random four-corner Cantor set. The lower bound of ${\rm Fav}({\cal R}_n)$ is $cn^{-1}$ again by Mattila \cite{m}. The interesting thing is Peres and Solomyak proved that, with high probability, the upper bound of ${\rm Fav}({\cal R}_n)$ is also constant
multiples of $n^{-1}$.  The random model of Peres and Solomyak convinced people that the $n^{-1}$ lower bound by Mattila is optimal in some general sense. Secondly, recall that the best known results for the non-random case are $n^{-1}\log n\lesssim {\rm Fav}({\cal K}_n)\lesssim n^{-\frac{1}{6}+\epsilon}$.  The difference between the decay rate of ${\rm Fav}({\cal K}_n)$ and ${\rm Fav}({\cal R}_n)$ inspires us to consider whether there are indeed some transition phenomenons between the deterministic graph and the random graph or not.  Such transition exists widely in ergodic theory and math-physics literatures, which is another motivation for our current paper and for some of the future projects.\\

Now we start to introduce the random Cantor set we want to study. Instead of the square, an equivalent model is to start with the unit disk. Namely, let ${\cal G}_0:=\{z\in\C:\ |z|\le 1\}$ be the unit disk. Replace ${\cal G}_0$ by
four subdisks of radius $1/4$ centered at $(0,\pm3/4),(\pm3/4,0)$, and iterate this self-similarly in each subdisk. More formally, let $T_j(z)=\frac{1}{4}z+\frac{3}{4}e^{{\rm i}\frac{\pi}{2} j},\ z\in\C,\ j=0,1,2,3$ and
\begin{equation}\label{Gn}
  {\cal G}_n=\bigcup_{j=0}^{3}T_j({\cal G}_{n-1}),\ \ n=1,2,3,\cdots.
\end{equation}
Then the quarter disk Cantor set $\mathcal{G}$ will be defined as $\mathcal{G}=\bigcap_{n=0}^{\infty}{\cal G}_n$. Clearly, the quarter disk Cantor set $\mathcal{G}$ is equivalent to the four-corner Cantor set $\mathcal{K}$. Almost all conclusions of ${\cal K}_n$ hold true for ${\cal G}_n$.

In this paper, we consider a random analogue of the above quarter disk model.
We start again with the unit disk ${\cal D}_0:={\cal G}_0$. In the first step, after placing the four subdisks with quarter radius centered at $(0,\pm3/4),(\pm3/4,0)$, we add further a random rotation to the picture, i.e., consider
\begin{equation}\label{D1}
{\cal D}_1(\omega_1)=\bigcup_{j=0}^{3}T^{\omega_1}_j(\calD_{0}),
\end{equation}
where $T^{\omega_1}_j(z)=\frac{1}{4}z+\frac{3}{4}e^{{\rm i}(\frac{\pi}{2} j-\omega_1)},\ z\in\C,\ j=0,1,2,3$ and $\omega_1$ is uniformly distributed random variable in $[0,\frac{\pi}{2}]$. There are four quarter subdisks in ${\cal D}_1(\omega_1)$. In the next step, we replace each of these four subdisks by four smaller disks randomly given by $\bigcup_{j=0}^{3}\frac{1}{4}T^{\omega_2}_j(\calD_{0})$ and denote the new figure (the union of $16$ subdisks with radius $1/16$) by ${\cal D}_2(\omega_1,\omega_2)$, where $T^{\omega_2}_j$ is defined in the same way as $T^{\omega_1}_j$ and $\omega_1,\omega_2$ are independent uniformly distributed random variables in $[0,\frac{\pi}{2}]$. Inductively, in the $k$-th  step, suppose we have ${\cal D}_{k-1}(\omega_1,\cdots,\omega_{k-1})$ which contains $4^{k-1}$ subdisks with radius $4^{-k+1}$, we replace each subdisk by the following group of 4 disks with radius $4^{-k}$:
\begin{equation}\label{Ck}
 {\cal C}_{k}( \omega_k)=\bigcup_{j=0}^{3}\frac{1}{4^{k-1}}T^{\omega_{k}}_j(\calD_{0}),
\end{equation}
where $T^{\omega_{k}}_j(z)=\frac{1}{4}z+\frac{3}{4}e^{{\rm i}(\frac{\pi}{2} j-\omega_{k})},\ z\in\C,\ j=0,1,2,3$ and $\omega_{k}$ is uniformly distributed random variable in $[0,\frac{\pi}{2}]$ and independent of $(\omega_1,\cdots,\omega_{k-1})$. The collection of all these $4^{k}$ subdisks of radius $4^{-k}$ is denoted as  ${\cal D}_{k}(\omega_1,\cdots,\omega_{k})\subset {\cal D}_{k-1}(\omega_1,\cdots,\omega_{k-1})$. The random analogue of ${\cal G}$ is now defined by ${\cal D}=\bigcap{\cal D}_k$ and ${\cal D}_k$ is a $4^{-k}$ neighborhood of ${\cal D}$.   Let us now focus on the $n$-th  step, let $\omega=(\omega_1,\omega_2,\cdots,\omega_n)\in[0,\frac{\pi}{2}]^n$ and denote $\calD_n(\omega)={\cal D}_{n}(\omega_1,\cdots,\omega_{n})$ for simplicity. Clearly, ${\cal D}_n(\omega)$ can be viewed as rotating each layer of ${\cal G}_n$ in (\ref{Gn}) by the angles $(\omega_1,\omega_2,\cdots,\omega_n)$ independently since ${\cal G}_n={\cal D}_n(0,\cdots,0)$. It is easy to check that the above rotation does not change the ${\cal H}^1$ measure\footnote{We refer readers more background about such self-similar Cantor sets, Hausdorff dimension/measure and orthogonal projection to \cite{f,m3} and references therein.}. Therefore, we have $\lim_n{\rm Fav}({\cal D}_n(\omega))=0$ for the same reason as for the deterministic case ${\cal G}_n$. Similarly, we are interested in a good estimate for ${\rm Fav}({\cal D}_n(\omega))\to 0$ as $n\to \infty$.

For any $\omega\in [0,\frac{\pi}{2}]^n$, the lower bound ${\rm Fav}(\calD_n(\omega))\ge c_{\omega} n^{-1}$ follows again from Mattila \cite{m}. Our main result is about the upper bound. Denote the expectation by $\E_{\omega}$, we have:
\begin{thm}\label{mainthm}
Suppose $\omega=(\omega_1,\omega_2,\cdots,\omega_n)$ are uniformly distributed i.i.d. random variables in $[0,\frac{\pi}{2}]$. Then there exists $C>0$ such that for any $\theta\in[0,\pi]$,
\begin{equation}\label{maineq1}
  \E_{\omega}\Big|{\rm Proj}_{\theta}\ \calD_n(\omega)\Big|\le \frac{C}{n}\ \ \  \textrm{for all}\ \  \  n\in\N.
\end{equation}
Consequently, we have
\begin{equation}\label{maineq2}
  \E_{\omega}\big[{\rm Fav}\big(\calD_n(\omega)\big)\big]\le \frac{C}{n}\ \ \  \textrm{for all}\ \  \  n\in\N
\end{equation}
and
\begin{equation}\label{maineq3}
  \liminf_{n}\ n\cdot{\rm Fav}\big(\calD_n(\omega)\big)<\infty\ \ \  \textrm{for all}\ \  \  n\in\N,\ \textrm{almost surely}.
\end{equation}
\end{thm}

\begin{rmk}\label{rmk:npv}
Combine the lower bound of Mattila with Theorem \ref{mainthm}, we see that the upper and lower bounds of ${\rm Fav}({\cal D}_n(\omega))$ are both $n^{-1}$ up to some constants in the average sense. In the deterministic case, the best known results are $n^{-1}\log n\lesssim {\rm Fav}({\cal G}_n)\lesssim n^{-\frac{1}{6}+\epsilon}$ proved in \cite{npv,bv}. It is still an open problem to
determine the exact rate of decay of $ {\rm Fav}({\cal G}_n)$.
\end{rmk}
\begin{rmk}\label{rmk:ps}
Our model is inspired by the random four-corner Cantor set ${\cal R}$ in \cite{ps} but in a slightly different taste. Peres and Solomyak considered a discrete distribution with 4 choices at each node in \cite{ps}. The uniformly distributed rotation we considered is more natural to the disk.   Our main results focus on uniform distribution for simplicity. The scheme potentially can be used to study many other i.i.d. distributions as long as the expectation of the overlap of two intervals can be bounded from below by square order (see Lemma \ref{IJ}). Secondly, our model is in some sense less random compared to ${\cal R}$. Recall in (\ref{Ck}), we replace all $4^{k-1}$ subdisks by the same ${\cal C}_k(\omega_k)$. We can actually replace each subdisk by different ${\cal C}_k(\omega_k^m),\ m=1,\cdots,4^{k-1}$ where $\omega_k^m$ are all i.i.d. random variables. This choice would resemble the randomness of ${\cal R}$ considered by Peres and Solomyak. As we saw in Remark \ref{rmk:npv}, the deterministic case is usually considered much harder than the random case. We do not intend to generalize our model to other i.i.d distribution with more randomness. The future goal is to study some less random model, i.e., only assume the sequence $\omega_k$ to have certain weaker/mild ergodicity rather than i.i.d.
\end{rmk}

The proof of the main theorem will be different from either \cite{npv} or \cite{ps}.
For the non-random case \cite{npv,lz,bv1,bv2,blv}, the estimates either rely on the tiling properties of the Cantor structures or the small value properties of certain trigonometric polynomials, see more details in \cite{npv,blv} and references therein. For the random case ${\cal R}_n$ in \cite{ps}, the proof relies on an argument involving
percolation on trees and some delicate inductions. In the current model, none of these are needed. The proof is simply based on estimating the average overlap of the projections and some inductive argument from the independency of the random variables. The independency of the random variables allow us to rotate the $k$-th layer by $\omega_k$ and take the expectation separately. What's more important, the rotations $(\omega_1,\cdots,\omega_n)$ are added from the 1st layer to the $n$-th layer one by one independently. Once we are done with the construction and fix the all inner rotations $(\omega_2,\cdots,\omega_n)$, the 1st/outer rotation  $\omega_1$ is actually a shift with respect to all the inner layers since it only moves the center of each sub image without changing the orientation/projection length.

Since $(\omega_1,\cdots,\omega_n)$ are i.i.d. uniform rotations, it is enough to prove (\ref{maineq1}) for $\theta=0$. (\ref{maineq2}) follows directly from (\ref{maineq1}) and the definition of the Favard length. Once we show (\ref{maineq2}), Fatou's lemma leads to (\ref{maineq3}). In the rest of the paper, we will denote ${\rm Proj}\, E={\rm Proj}_{\theta=0}\, E$ for short. In section \ref{keysec}, we introduce the key inductive Lemma \ref{keylem} and  prove (\ref{maineq1}). In section \ref{pf}, we will prove the key inductive Lemma \ref{keylem}. In the last section, we discuss briefly the generalization of the quarter disk Cantor set $\calD$ to $1/d$ disk Cantor set $\calD^d$, which is constructed with $d$ subdisks of radius $1/d$ for any $d\ge3$. We will see in Section \ref{secDd} that all the conclusions and  arguments for $\calD$ hold true for $\calD^d$.\\

\noindent {\bf Acknowledgment.} I wish to thank Alexander Volberg for suggesting to me the problem
of studying the decay rate of the Favard length and the random Cantor disk. I appreciate his support from the very beginning of this project and for many useful suggestions and comments. I would also like to thank Fedor Nazarov for his support and hospitality at Kent State University where the work was done. I am grateful to his help and important suggestions concerning the key Lemma \ref{keylem}.
\section{Key inductive lemma and Proof of the main result}\label{keysec}
Let ${\cal G}_n$ be the unrotated graph given as in (\ref{Gn}). Now let's rotate ${\cal G}_n$ layer by layer randomly to get ${\cal D}_n$. We consider the expectation of the projection length from the $n$-th layer back to the first layer.
One simple observation is: there will be more overlap for the projection in the average sense under rotation. More importantly, the overlap is of square order of the original projection length. The quantitative estimate on the average overlap leads to the $1/n$ decay inductively. For $k=1,2,\cdots,n$, denote
\begin{equation}\label{omegak}
 \bar{\omega}^k=(\omega_{n-k+1},\cdots,\omega_{n-1},\omega_n)\in [0,\frac{\pi}{2}]^{k},
\end{equation}
and consider the expected projection length onto the horizontal axis:
\begin{equation}\label{Ek}
E_k=\E_{ \bar{\omega}^k} \Big|{\rm Proj}\ {\cal D}_k( \bar{\omega}^k)\Big|,
\end{equation}
where ${\cal D}_k( \bar{\omega}^k)$ is constructed step by step as in (\ref{D1}), (\ref{Ck}) from $\omega_{n-k+1}$ to $\omega_n$.\\

Our key inductive lemma is as follows:
\begin{lemma}\label{keylem}
Suppose $ {\omega}=(\omega_{1},\cdots,\omega_n)$ are unifromly distributed i.i.d random variables in $[0,\frac{\pi}{2}] $, then there exists absolute constant $0<c<\frac{1}{4}$ such that for any $2\le n\in\N$,
\begin{equation}\label{keyineq}
  E_k\le E_{k-1}-cE^2_{k-1}, \ \ k=2,\cdots,n.
\end{equation}
\end{lemma}

\noindent \textbf{Proof of Theorem \ref{mainthm}:} Let $c$ be given in Lemma \ref{keylem} and $C=c^{-1}$. According to the definition, $$E_1=\E_{ \omega_n} \Big|{\rm Proj}\calD_{1}( \omega_n)\Big|.$$
Clearly, $E_1\le2<C$.
 Now assume that $E_{k-1}<\frac{C}{k-1}$ and do induction from $k=2$ to $k=n$. (\ref{keyineq}) gives immediately that
 \begin{equation}\label{Ekdecay}
   E_{k}<\frac{C}{k-1}-\frac{C}{(k-1)^2}=C\frac{k-2}{(k-1)^2}<\frac{C}{k}.
 \end{equation}
 Therefore, (\ref{Ekdecay}) holds for all $k=1,\cdots,n$. In particular,
$$\E_\omega \Big|{\rm Proj} \, \calD_n(\omega)\Big|=E_n<\frac{C}{n}.$$ \qed

\begin{rmk}
Lemma \ref{keylem} holds true for projections with any angle $\theta$. The above proof for (\ref{maineq1}) will be exactly the same for any other $\theta$. See more discussion in the end of next section.
\end{rmk}

\begin{rmk}
Once ${\cal D}_n(\omega_1,\cdots,\omega_n)$ is constructed, we want to study the structure layer by layer. If we look at the $(n-k)$-th layer, there are actually $4^{n-k}$ copies of the figure
$$4^{-(n-k)}{\cal D}_k( \omega_{n-k+1},\cdots,\omega_n).$$
All these copies only differ by a shift of the centers and have exactly the same orientation and the projection length. That's why we can actually drop the scaling factor $4^{k-n}$ and consider $E_k$ only.
\end{rmk}

\section{Proof of the key inductive Lemma \ref{keylem}}\label{pf}
Let $\bar\omega^k$ be given as in (\ref{omegak}). Denote
\begin{equation}\label{Lk}
  L_k( \bar{\omega}^k)=\Big|{\rm Proj}\, {\cal D}_{k}( \bar{\omega}^k)\Big|
\end{equation}
Write $\bar\omega^k=(\tilde \omega,\omega')$ for short, where $\tilde \omega=\omega_{n-k+1}$, $\omega'=(\omega_{n-k+2},\cdots,\omega_n)$. We have $\E_{\bar\omega^k}L_k=\E_{\tilde \omega}\big(\E_{\omega'}L_k\big)$ since $\bar\omega^k$ are i.i.d. random variables.  Now let us fix $\omega'$ and consider the expectation on $\tilde \omega$ first. The key to prove (\ref{keyineq}) is the following lemma on $\E_{\tilde \omega}L_k(\tilde \omega,\omega')$:
\begin{lemma}\label{omegaklem}
Assume $\tilde \omega$ obeys uniform distribution on $[0,\frac{\pi}{2}]$. There is an absolute constant $0<c<\frac{1}{4}$ such that for any  $k=2,\cdots,n$, and $\omega'\in[0,\frac{\pi}{2}]^{k-1}$, we have
\begin{equation}\label{omegakeq}
  \E_{\tilde \omega}L_k(\tilde \omega,\omega')\le L_{k-1}(\omega')-cL^2_{k-1}(\omega').
\end{equation}
\end{lemma}

Once (\ref{omegakeq}) is established, integrate both sides with respect to $\omega'$ to obtain
 \begin{eqnarray}
 \E_{(\tilde \omega,\omega')}L_k(\tilde \omega,\omega')&\le&  \E_{\omega'}L_{k-1}(\omega')-c\E_{\omega'}L^2_{k-1}(\omega')\nonumber\\
 &\le& \E_{\omega'}L_{k-1}(\omega')-c\big(\E_{\omega'}L_{k-1}\big)^2, \label{expLk}
 \end{eqnarray}
where the last line follows from applying Cauchy-Schwarz to $\big(\E_{\omega'}L_{k-1}\big)^2\le \E_{\omega'}L^2_{k-1}$.
Notice that $\omega'=\bar\omega^{k-1}$ and $E_{k-1}=\E_{\bar\omega'}L_{k-1}(\omega')$, (\ref{expLk}) immediately implies (\ref{keyineq}). In the rest of section, we focus on proving (\ref{omegakeq}).

Given $(\tilde \omega,\omega')$, let's look at the structure of ${\cal D}_k(\bar\omega^k)={\cal D}_k(\tilde \omega,\omega')$.  ${\cal D}_k(\tilde \omega,\omega')$ contains four groups of images, each group is a copy of $\frac{1}{4}{\cal D}_{k-1}(\omega')$. Let's call these copies $U^1,U^2,U^3,U^4$ (in clockwise order).
Fix $\omega'$, for any $\tilde \omega$, all these four $U^j(\tilde\omega,\omega')$ have the same projection length $\frac{1}{4}|{\rm Proj}\,{\cal D}_{k-1}(\omega')|=\frac{1}{4}L_{k-1}(\omega')$. Fix $\omega'$ and only rotate
\footnote{As we have explained in Section \ref{sec:intro}, it is actually  a ``shift'' of the four $U^j$. When we say ``rotate'' the outside layer, we only map the center of each $U^j(0,\omega')$ by $z\mapsto\frac{3}{4}e^{-{\rm i}\tilde\omega}z$ and then move the entire $U^j$ without changing its inner orientation given by $\omega'$. Most importantly, such ``rotation/shift'' does not change the projection length of each $U^j$.} the outside layer of ${\cal D}_k(\tilde \omega,\omega')$ by $\tilde \omega$. It's not hard to believe that the expected projection length of the rotated graph will decrease since there will be more overlap in the average sense. More precisely, for any $\tilde\omega$, any point on the horizontal axis can be covered by projections of at most two $U^j$, therefore,
\begin{eqnarray*}
\big|{\rm Proj}\,{\cal D}_k(\bar\omega^k)\big|&\le& \big|\bigcup_{j=1}^4{\rm Proj}\, U^j\big|-\sum_{i,j=1,i\neq j}^4\big|{\rm Proj}\,U^i\cap {\rm Proj}\,U^j\big|\\
&\le& \sum_{j=1}^4\big|{\rm Proj}\, U^j(\tilde \omega,\omega')\big|-\big|{\rm Proj}\,U^1(\tilde \omega,\omega')\cap {\rm Proj}\,U^2(\tilde \omega,\omega')\big| \\
&=& \sum_{j=1}^4\frac{1}{4}|{\rm Proj}\,{\cal D}_{k-1}(\omega')|-\big|{\rm Proj}\,U^1(\tilde \omega,\omega')\cap {\rm Proj}\,U^2(\tilde \omega,\omega')\big|,
\end{eqnarray*}
i.e.,
$$L_k(\tilde \omega,\omega')\le L_{k-1}(\omega')-\big|{\rm Proj}\,U^1(\tilde \omega,\omega')\cap {\rm Proj}\,U^2(\tilde \omega,\omega')\big|.$$
Take expectation with respect to $\tilde \omega$, we have
\begin{equation}\label{EkEk+1}
  \E_{\tilde\omega}L_k(\tilde\omega,\omega')\le L_{k-1}(\omega')-\E_{\tilde\omega}\Big(\big|{\rm Proj}\,U^1(\tilde \omega,\omega')\cap {\rm Proj}\,U^2(\tilde \omega,\omega')\big|\Big)
\end{equation}
To prove (\ref{omegakeq}), it is enough to bound the average overlap of the projections from below. Suppose ${\cal D}_n(0,\omega')$ is the unrotated graph. Without loss the generality,  assume $U^1(0,\omega')$ is the copy of ${\cal D}_{k-1}(\omega')$ in the North and $U^2(0,\omega')$ is the copy in the East. Denote $I={\rm Proj}\,U^1(0,\omega')$ and $J={\rm Proj}\,U^2(0,\omega')=I+\frac{3}{4}$. Now ``rotate/shift'' the centers of $U^1(0,\omega')$ and $U^2(0,\omega')$ by $z\mapsto\frac{3}{4}e^{-{\rm i}\tilde\omega}z$. The new projections now become  $I(\tilde \omega):={\rm Proj}\,U^1(\tilde\omega,\omega')=I+\frac{3}{4}\sin\tilde \omega$ and  $J(\tilde \omega):={\rm Proj}\,U^2(\tilde\omega,\omega')=I+\frac{3}{4}\cos\tilde \omega$.
The lower bound of $\E_{\tilde\omega}\Big(\big|{\rm Proj}\,U^1(\tilde \omega,\omega')\cap {\rm Proj}\,U^2(\tilde \omega,\omega')\big|\Big)$ follows from the following technical lemma:
\begin{lemma}\label{IJ}
For any $a>0$ and $I\subset [-a,a]$
\begin{equation}\label{IJineq}
\int_{0}^{\pi/2}\Big|(I+3a\sin\theta) \cap (I+3a\cos\theta)\Big|{\rm d}\theta \ge \frac{1}{6\sqrt2a} |I|^2
\end{equation}
\end{lemma}

\noindent \textbf{Proof of the Lemma:}  As $\theta$ increases from $0$ to $\pi/2$, $I+3a\sin\theta$ moves to the right and $I+3a\cos\theta$ moves to the left.  Let $\theta_\ast\in[0,\frac{\pi}{4}]$ be such that $a+3a\sin\theta_\ast=-a+3a\cos\theta_\ast$ where the two sets start to overlap. It is easy to check that $(I+3a\sin\theta)\cap(I+3a\cos\theta)=\emptyset$ for $\theta\in [0,\theta_\ast)\cup(\frac{\pi}{2}-\theta_\ast,\frac{\pi}{2}]$ and the two intervals coincide at $\theta=\frac{\pi}{4}$.
And $(I+3a\sin\theta)\cap(I+3a\cos\theta)\subset [-a+3a\cos\theta,a+3a\sin\theta]$ for $\theta\in[\theta_\ast,\frac{\pi}{2}-\theta_\ast]$. Therefore, by symmetry, it is enough to estimate
\begin{eqnarray}
 &&\int_{\theta_\ast}^{\pi/4}\int\chi_{I+3a\sin\theta}(x)\,\chi_{I+3a\cos\theta}(x)\ {\rm d}x\, {\rm d}\theta \nonumber \\
&=&\int_{\theta_\ast}^{\pi/4}\int_{-a+3a\cos\theta}^{a+3a\sin\theta}\chi_{I}(x-3a\sin\theta)\,\chi_{I}(x-3a\cos\theta)\ {\rm d}x\, {\rm d}\theta .\label{uv}
\end{eqnarray}
By changing variables $u=x-3a\sin\theta$ and $v=x-3a\cos\theta$,  we have

$$\theta_\ast\le \theta\le \pi/4,\ -a+3a\cos\theta\le x\le a+3a\sin\theta
\Longrightarrow -a\le v\le u\le a,\ \ {\rm so}$$
\begin{equation}\label{uv2}
  (\ref{uv}) =\iint_{-a\le v\le u\le a}\chi_{I}(u)\chi_{I}(v)\
  \Big|\frac{\partial(\theta,x)}{\partial(u,v)}\Big|{\rm d}u {\rm d}v.
\end{equation}
Direct computation shows
$$3a\le \Big|\frac{\partial(u,v)}{\partial(\theta,x)}\Big|=3a\sqrt 2\sin(\theta+\frac{\pi}{4})\le 3a\sqrt2.$$
Therefore,
$$(\ref{uv2})\ge \frac{1}{3a\sqrt 2} \iint_{-a\le v\le u\le a}\chi_{I}(u)\chi_{I}(v)\
 {\rm d}u {\rm d}v= \frac{1}{3\sqrt 2 a}\cdot \frac{1}{2}|I|^2
$$

 The last equality comes from the following symmetry of the integral:
 $$\int_{-a}^{a}\int_{v}^{a}
  \chi_{I}(u)\ \chi_{I}(v)\ {\rm d}u\,  {\rm d}v
 =\int_{-a}^{a} {\rm d}u\int_{-a}^{u}\chi_{I}(u)\chi_{I} (v)\ {\rm d}v=\int_{-a}^{a} \int_{-a}^{v}\chi_{I}(u)\chi_{I} (v)\ {\rm d}u\,{\rm d}v.
 $$
Sum the left and right integrals up, to get
 $$\int_{-a}^{a}\int_{v}^{a}
  \chi_{I}(u)\ \chi_{I}(v)\ {\rm d}u\,  {\rm d}v
 =\frac{1}{2}\int_{-a}^{a}\int_{-a}^{a}\chi_{I}(u)\chi_{I} (v)\ {\rm d}u\ {\rm d}v
 =\frac{1}{2}|I|^2.$$ \qed
 
 \begin{rmk}
 This is the case for the projection onto the horizontal axis. It is easy to see that when we consider the general projection with angle $\theta_0$, i.e., $ L^{\theta_0}_k( \bar{\omega}^k)={\rm Proj}_{\theta_0} \, {\cal D}_k( \bar{\omega}^k)$. The only change in Lemma \ref{IJ} is to estimate the integral over $[\theta_0,\theta_0+\frac{\pi}{2}]$. The lower bound will the same since the two adjacent $U^j$ differ by an angle $\pi/2$.
 \end{rmk}
 Now apply the above lemma to $I(\tilde \omega)$ and $J(\tilde \omega)$ with $a=1/4$.
 Notice that $|I(\omega')|=\frac{1}{4}L_{k-1}(\omega')$, so (\ref{IJineq}) implies that
 \begin{equation}\label{IJLk-1}
 \int_{0}^{\pi/2}\big|I(\tilde \omega)\cap J(\tilde \omega)\big|\ {\rm d}\tilde \omega
 \ge \frac{\sqrt2}{3}|I(\omega')|^2= \frac{\sqrt2}{48} L^2_{k-1}(\omega').
 \end{equation}

Clearly, (\ref{omegakeq}) follows from (\ref{EkEk+1}) and (\ref{IJLk-1}). This completes the proof of Lemma \ref{omegaklem}. \qed

\section{General random disk model ${\mathcal D}^d(\omega)$ with ratio $1/d$, }\label{secDd}
An advantage for considering the disk rather than regular polygons is that it is  easier to be generalized to a self-similar set with higher degree. In \cite{bv2}, Bond and Volberg considered the general disk model ${\cal G}^d$ with ratio $d^{-1}$ for all $3\le d\in\N$,
\begin{equation}\label{Gd}
  {\cal G}^d_n=\bigcup_{j=0}^{d-1}S_j({\cal G}^d_{n-1}),\ \ n=1,2,3,\cdots,
\end{equation}
where ${\cal G}^d_0$ is the unit disk, $S_j(z)=\frac{1}{d}z+\frac{d-1}{d}e^{{\rm i}\frac{2\pi}{d} j},\ z\in\C,\ j=0,1,2,d-1$ and ${\cal G}^d=\bigcap_{n}{\cal G}^d_n$. Following the arguments developed in \cite{npv}, they can show that ${\rm Fav}({\cal G}^3_n)$ has similar power decay as ${\cal K}_n$. For $d\ge5$, no power decay was found; the best result is ${\rm Fav}({\cal G}^d_n)\le C e^{-c\sqrt{\log n}}$ in \cite{bv2}.

The random analogue of ${\cal G}^d$, namely, $\calD^d(\omega_1,\omega_2,\cdots)$ can be defined exactly in the same way as in (\ref{D1}), (\ref{Ck}) for $\omega_k\in[0,\frac{2\pi}{d}],k=1,2,\cdots$. For uniformly distributed random variables $\omega=(\omega_1,\omega_2,\cdots,\omega_n)\in[0,\frac{2\pi}{d}]^n$, suppose we have already constructed ${\cal D}^d_n(\omega)$ through the n-step random rotations $\omega$.  Using the same notation as in (\ref{Ek}) and (\ref{Lk}), we can now apply Lemma \ref{IJ} with $a=1/d$ to the projection length of ${\cal D}^d_n(\omega)$. Absolutely the same reasoning as for $d=4$ proves
similar results as in Lemma \ref{omegaklem} and \ref{keylem} for any $d$:
\begin{lemma}
Let $\tilde \omega=\omega_{n-k+1}$ and $\omega'=(\omega_{n-k+2},\cdots,\omega_n)$. There exists absolute constant $0<c_d<1/d$ such that for any $2\le n\in\N$, $k=2,\cdots,n$ and $\omega'\in[0,\pi/2]^{k-1}$, we have
\begin{equation}
\E_{\tilde \omega}\big|{\rm Proj}\, {\cal D}^d_k(\tilde \omega,\omega')\big|\le
\big|{\rm Proj}\, {\cal D}^d_{k-1}(\omega')\big|-c_d\big|{\rm Proj}\, {\cal D}^d_{k-1}(\omega')\big|^2
\end{equation}
and
\begin{equation}
\E_{(\tilde\omega,\omega')}\big|{\rm Proj}\, {\cal D}^d_k(\tilde \omega,\omega')\big|\le
\E_{\omega'}\big|{\rm Proj}\, {\cal D}^d_{k-1}(\omega')\big|
-c_d\Big[\E_{\omega'}\big|{\rm Proj}\, {\cal D}^d_{k-1}(\omega')\big|\Big]^2.
\end{equation}
\end{lemma}

With the above inductive lemma, clearly, we have the same decay for ${\rm Fav}(\calD^d_n)$ as in Theorem \ref{mainthm}.
\begin{thm}\label{Ddthm}
Suppose $d\ge3, d\in\N$ and $\omega=(\omega_1,\omega_2,\cdots,\omega_n)$ are uniformly distributed i.i.d. random variables in $[0,\frac{2\pi}{d}]$. Then there exists $C_d>0$ such that for any $\theta\in[0,\pi]$,
\begin{equation}
  \E_{\omega}\Big|{\rm Proj}_{\theta}\ \calD_n(\omega)\Big|\le \frac{C_d}{n}\ \ \  \textrm{for all}\ \  \  n\in\N.
\end{equation}
Consequently, we have
\begin{equation}
  \E_{\omega}\big[{\rm Fav}\big(\calD^d_n(\omega)\big)\big]\le \frac{C_d}{n}\ \ \  \textrm{for all}\ \  \  n\in\N
\end{equation}
and
\begin{equation}
  \liminf_{n}\ n\cdot{\rm Fav}\big(\calD^d_n(\omega)\big)<\infty\ \ \  \textrm{for all}\ \  \  n\in\N,\ \textrm{almost surely}.
\end{equation}
\end{thm}

\begin{rmk}
It is interesting to see that for $d\ge5$, there is even no power decay for the deterministic case. While in the random case, we can easily get $1/n$ decay for any $d$  with no essential difference in the proof.
\end{rmk}


{}

\vspace{1cm}
Shiwen Zhang, 

Dept. of Math., Michigan State University. 

E-mail address: zhangshiwen@math.msu.edu

\end{document}